\documentclass[12]{article}

\usepackage[T2A]{fontenc}
\usepackage[cp1251]{inputenc}
\usepackage[russian,english]{babel}
\usepackage[tbtags]{amsmath}
\usepackage{amssymb,amsfonts}
\usepackage{mathrsfs}
\usepackage[12pt]{extsizes}

\begin{document}

\begin{center}
{\large{ON TWO-DIMENSIONAL DIRICHLET SPECTRUM}}

 Renat Akhunzhanov

 Denis Shatskov

 {\textit{Astrakhan State University}}
\end{center}

\begin{abstract}
{\textit{We define two-dimensional Dirichlet spectrum (with
respect to Euclidean norm) as
$$
\mathbb{D}_2=\{\lambda\in\mathbb{R}\mid\exists\mathbf{v}=(v_1,v_2)\in\mathbb{R}^2:\;
\limsup\limits_{t\rightarrow\infty}
{t\cdot\psi_{\mathbf{v}}^2(t)}=\lambda\},
$$
where
$$
\psi_{\mathbf{v}}(t)=\min\limits_{1\leqslant q\leqslant
t}\sqrt{\|q v_1\|^2+\|q v_2\|^2}
$$
is the
two-dimensional ``irrationality measure function''.\\
Our main result  states the equality
$$
\mathbb{D}_2=\left[0,\frac{2}{\sqrt{3}}\right].
$$}}
\end{abstract}

\vspace{10ex}

First of all we recall well known one-dimensional facts. Let
$\alpha\in \mathbb{R}$ be an irrational number.

Irrationality measure function for $\alpha$ is defined as
$$
\psi_\alpha(t)=\min\limits_{1\leqslant q\leqslant t}\|q\alpha\|.
$$

This function can be also described in terms of continued
fractions.

Consider continued fraction expansion $\alpha=[a_0;a_1,a_2,\dots]$
and convergents $\frac{p_n}{q_n}=[a_0;a_1,a_2,\dots,a_n]$.

Then $\psi_\alpha(t)$ is  piecewise constant decreasing function
satisfying
$$
\psi_\alpha(t)=\|q_n\alpha\|, \mbox{ for } q_n\leqslant t<q_{n+1}.
$$

Dirichlet spectrum is defined as
$$
\mathbb{D}=\{\lambda\in\mathbb{R}\mid\exists\alpha\in\mathbb{R}:\;
\limsup\limits_{t\rightarrow\infty} {t\psi_\alpha(t)}=\lambda\}.
$$
An equivalent definition is as follows:
$$
\mathbb{D}=\{\lambda\in\mathbb{R}\mid\exists\alpha\in\mathbb{R}:\;
\limsup\limits_{n\rightarrow\infty}
{q_{n+1}\|q_n\alpha\|}=\lambda\}.
$$
We note that the expression under $\limsup$ has a clear  geometric
meaning.
 We consider the parallelogram
$$
\Pi(\alpha, X, R)= \left\{(x,y)\in\mathbb{R}^{2}\mid x\in X,\;
|x\alpha - y|\leqslant R\right\}.
$$
Put $\Pi_n=\Pi(\alpha, [0, q_n], R_n)$,  and let $S_n=S(\Pi_n)$ be
the
  area of the $\Pi_n$.
So
$$q_{n+1}\|q_n\alpha\|=S_{n+1}=S(\Pi_{n+1}).$$

From Minkowski convex body theorem one can easily see that
$$
\mathbb{D}\subset[0,1]
$$
G. Szekeres (\cite{200}, 1937) showed that
$$
\mathbb{D}\subset\left[\frac{5+\sqrt{5}}{10},1\right]=[0.72\dots,1].
$$

It is known that $\frac{5+\sqrt{5}}{10}$ is the smallest point in
$ \mathbb{D}$. The discrete part of $ \mathbb{D} $ was studied by
B. Divis (\cite{Di},1972). He also proved that  there exists
$d\in(0,1)$, such that $\mathbb{D}\supset[d,1]$. Define
$d^*=\inf\{d\mid\mathbb{D}\supset[d,1]\}.$ V.A. Ivanov
(\cite{100},1980) proved that
$$
d^*\in\left[\frac{3\sqrt{5}-5}{2},\frac{38+6\sqrt{2}}{49}\right]=[0.85\dots,0.94\dots]
$$
In  1978 V.A. Ivanov (\cite{100}, 1980) also showed that
$$
\mbox{mes}\left(\mathbb{D}\cap\left[0,\frac{4+3\sqrt{3}}{11}\right)\right)=0
$$
Note that
$$
\frac{4+3\sqrt{3}}{11} = 0.83\dots
$$

The main tool for all the results behind is continued fraction
representation of $\alpha$.
  The  basic relation is the   equality
$$
q_{n+1}\cdot\|q_n\alpha\|=\frac{1}{1+\frac{1}{\alpha_{n+2}\alpha_{n+1}^{**}}},
$$
where $\alpha_n=[a_n;a_{n+1},a_{n+2},\dots]$ and
$\alpha_n^{**}=[a_n;a_{n-1},a_{n-2},\dots,a_1]$.

Up to our knowledge, the complete structure of one-dimensional
Dirichlet spectrum is not clear.

\vskip+0.4cm

Our paper is devoted to two-dimensional simultaneous Diophantine
approximations with respect to the Euclidean norm. Of course in
the two-dimensional case there is no such a tool as contunued
fractions. However we are able to  define the structure of the
dwo-dimensional Dirichlet spectrum completely.

For a vector $\mathbf{v}=(v_1, v_2)\in\mathbb{R}^2$ such that
 $1, v_1, v_2$ are linearly independent over $\mathbb{Z}$
we define two-dimensional ``irrationality measure function``
$$
\psi_{\mathbf{v}}(t)=\min\limits_{1\leqslant q\leqslant
t}\sqrt{\|q v_1\|^2+\|q v_2\|^2}.
$$
This function has jumps in points known as ''best
approximations''. They are defined inductively by
\\
 $q_0=1$,\\
$q_{n+1}=\min\left\{q\in\mathbb{N}\mid
\sqrt{\|q v_1\|^2+\|q v_2\|^2}<\sqrt{\|q_n v_1\|^2+\|q_n v_2\|^2}\right\}$,\\
$\mathbf{p}_n=(p_{n,1},p_{n,2}): \|q_n v_1\|=|q_n v_1-p_{n,1}|,
\|q_n v_2\|=|q_n v_2-p_{n,2}|$

We recall an alternative definition. Let the vector
${\bf{v}}\in\mathbb{R}^2$, the set $X\subset\mathbb{R}$ and $R>0$
be fixed. We denote
$$
\Pi({\bf{v}}, X, R)=
\left\{(q,{\bf{p}})\in\mathbb{R}\times\mathbb{R}^{2}\mid q\in X,\;
|q{\bf{v}} - {\bf{p}}|\leqslant R\right\}.
$$
  For a given vector $\mathbf{v}\in\mathbb{R}^2$
and the sequence
$$
\mathcal{Z}:\mathbf{w}_n=(q_n,\mathbf{p}_n)\in\mathbb{N}\times\mathbb{Z}^{2},\qquad
n\geqslant0
$$
we denote $R_0=1$, $R_n=|q_{n-1}\mathbf{v}-\mathbf{p}_{n-1}|$,
$\Pi_n=\Pi(\mathbf{v}, [0, q_n], R_n)$, $V_n=vol(\Pi_n)$. The
sequence $\mathcal{Z}$ will be   the sequence of the best
approximations for $\mathbf{v}\in\mathbb{R}^2$ with respect to
Euclidean norm, if \\
1) $q_0=1$\\
2) $\left(\mbox{int}\Pi_n\right)\cap\mathbb{Z}^3=\O$\\
3) $q_{n+1}>q_n$\\
4) $R_{n+1}<R_n$.

We define two-dimensional Dirichlet spectrum (with respect to
Euclidean norm) as
$$
\mathbb{D}_2=\{\lambda\in\mathbb{R}\mid\exists\mathbf{v}\in\mathbb{R}^2:\;
\limsup\limits_{t\rightarrow\infty}
{t\cdot\psi_{\mathbf{v}}^2(t)}=\lambda\}.
$$
Of course an equivalent definition will be
$$
\mathbb{D}_2=\{\lambda\in\mathbb{R}\mid\exists\mathbf{v}\in\mathbb{R}^2:\;
\limsup\limits_{n\rightarrow\infty}
{q_{n+1}\cdot|q_n\mathbf{v}-\mathbf{p}_n|^2}=\lambda\}.
$$
This definition has a clear geometric meaning as
$q_{n+1}\cdot|q_n\mathbf{v}-\mathbf{p}_n|^2=\frac{1}{\pi}V_{n+1}=\frac{1}{\pi}vol(\Pi_{n+1})$.

Minkowski convex body theorem trivially leads to
$$
\mathbb{D}_2\subset\left[0,\frac{4}{\pi}\right]
$$
K. Mahler's theorem on the critical determinant of a
three-dimensional cylinder implies
$$
\mathbb{D}_2\subset\left[0,\frac{2}{\sqrt{3}}\right]
$$

The main result of the present paper states that
$$
\mathbb{D}_2=\left[0,\frac{2}{\sqrt{3}}\right]
$$

In fact we prove more general result.

{\bf Theorem.}\\
Consider an arbitrary sequence
 $\{\Delta_n\}_1^{\infty}$ of open intervals
from the segment
 $\left[0,\frac{2}{\sqrt{3}}\right]$. Then there exists an uncountable family of vectors
 $\mathbf{v}\in\mathbb{R}^2$ such that
$$
q_{n+1}\cdot|q_n\mathbf{v}-\mathbf{p}_n|^2\in\Delta_n\;\;\;\forall
n\in\mathbb{N}
$$

The rest of the paper is written in Russian. We have submitted an
English version  to the Moscow Journal of Combinatorics and Number
Theory.

\newpage
\selectlanguage{russian}

\begin{center}
{\bf{Р.К. Ахунжанов, Д.О. Шацков}}\\

{\Large\bf{О двумерном спектре Дирихле}}
\end{center}

Пусть $\alpha\in \mathbb{R}$ иррациональное число. Определим
``функцию меры иррациональности'' числа~$\alpha$:
$$
\psi_\alpha(t)=\min\limits_{1\leqslant q\leqslant t}\|q\alpha\|.
$$
Рассмотрим разложение числа $\alpha$ в цепную дробь
$$
\alpha=[a_0;a_1,a_2,\dots].
$$
Пусть $\frac{p_n}{q_n}=[a_0;a_1,a_2,\dots,a_n]$ подходящая дробь
для $\alpha$.

Функция $\psi_\alpha(t)$ есть кусочно постоянная убывающая
функция:
$$
\psi_\alpha(t)=\|q_n\alpha\|, \mbox{ при } q_n\leqslant t<q_{n+1}
$$

Определим Спектр Дирихле (одномерный случай)
$$
\mathbb{D}=\{\lambda\in\mathbb{R}\mid\exists\alpha\in\mathbb{R}:\;
\limsup\limits_{t\rightarrow\infty}
{t\cdot\psi_\alpha(t)}=\lambda\}.
$$
Непосредственно из определения спектра Дирихле и свойств ``функции
меры иррациональности'' числа $\alpha$ получаем:
$$
\mathbb{D}=\{\lambda\in\mathbb{R}\mid\exists\alpha\in\mathbb{R}:\;
\limsup\limits_{n\rightarrow\infty}
{q_{n+1}\cdot\|q_n\alpha\|}=\lambda\}.
$$
Заметим, что $q_{n+1}\cdot\|q_n\alpha\|=\frac{1}{2}S(\Pi)$, где
$S(\Pi)$ это площадь параллелограмма

$$ \Pi=\{(x,y)\in\mathbb{R}^2\mid
 x\in[1, q_{n+1}], |x\alpha-y|\leqslant\|q_n\alpha\|\}
$$

Из теоремы Минковского о выпуклом теле следует тривиальный
результат:
$$
\mathbb{D}\cap(1,\infty)=\varnothing.
$$

Нетривиальные результаты:
$$
\mathbb{D}\subset\left[\frac{5+\sqrt{5}}{10},1\right]\;\;\;\mbox{
G. Szekeres (1937)}
$$

$$
\mbox{mes}\left(\mathbb{D}\cap\left[0,\frac{4+3\sqrt{3}}{11}\right)\right)=0\;\;\;\mbox{
В. А. Иванов (1978)}
$$

$$
d^*\in\left[\frac{3\sqrt{5}-5}{2},\frac{38+6\sqrt{2}}{49}\right],\;\;\;\mbox{
где } d^*=\inf\{d\mid\mathbb{D}\supset[d,1]\} \;\;\;\mbox{ В. А.
Иванов (1980)}
$$

Все вышеприведенные нетривиальные результаты были получены при
помощи аппарата цепных дробей и основаны на применении следующей
формулы:
$$
q_{n+1}\cdot\|q_n\alpha\|=\frac{1}{1+\frac{1}{\alpha_{n+2}\alpha_{n+1}^{**}}},
$$
где $\alpha_n=[a_n;a_{n+1},a_{n+2},\dots]$ и
$\alpha_n^{**}=[a_n;a_{n-1},a_{n-2},\dots,a_1]$.

{\bf Замечание.} В данной работе, для простоты изложения, мы даем
все определения и соответствующие результаты в евклидовой норме. По
всей видимости, аналогичные результаты можно получить и для других
норм.

{\bf Определение.} Для вектора $\mathbf{v}=(v_1,
v_2)\in\mathbb{R}^2$ мы определим двумерную ``функцию меры
иррациональности'', в евклидовой норме, следующим образом
$$
\psi_{\mathbf{v}}(t)=\min\limits_{1\leqslant q\leqslant t}\sqrt{\|q
v_1\|^2+\|q v_2\|^2}
$$

{\bf Определение.} Двумерным спектром Дирихле, в евклидовой норме,
называется множество
$$
\mathbb{D}_2=\left\{\lambda\in\mathbb{R}\mid\exists\mathbf{v}\in\mathbb{R}^2:\;
\limsup\limits_{t\rightarrow\infty}
{t\cdot\psi_{\mathbf{v}}^2(t)}=\lambda\right\}.
$$

Забегая вперед, можно сказать, что теорема Минковского о выпуклом
теле дает тривиальный результат
$$
\mathbb{D}_2\subset\left[0,\frac{4}{\pi}\right],
$$
а теорема Малера о критическом определителе трехмерного цилиндра
дает следующий результат:
$$
\mathbb{D}_2\subset\left[0,\frac{2}{\sqrt{3}}\right]
$$

В данной работе доказан следующий результат

{\bf Теорема 1.}\\
$$
\mathbb{D}_2=\left[0,\frac{2}{\sqrt{3}}\right].
$$

{\bf Определение.} Последовательность векторов
$$
\mathcal{Z}:\mathbf{w}_n=(q_n,\mathbf{p}_n)\in\mathbb{Z}^{3},\qquad
n\geqslant0
$$
называется последовательностью наилучших приближений для вектора
$\mathbf{v}\in\mathbb{R}^2$, в евклидовой норме, если\\
1) $q_0=1$;\\
2) $q_{n+1}=\min\left\{q\in\mathbb{N}\mid
\sqrt{\|q v_1\|^2+\|q v_2\|^2}<\sqrt{\|q_n v_1\|^2+\|q_n v_2\|^2}\right\}$ $(n\geqslant 0)$;\\
3) $\mathbf{p}_n=(p_{n,1},p_{n,2}): \|q_n v_1\|=|q_n v_1-p_{n,1}|,
\|q_n v_2\|=|q_n v_2-p_{n,2}|$ $(n\geqslant 0)$.

{\bf Замечание.} В общем случае, последовательность наилучших
приближений не всегда определена однозначно. Более точно,
последовательность $\{q_n\}_{n=0}^\infty$ определена однозначно, а
вот последовательность $\{\mathbf{p}_n\}_{n=0}^\infty$ может быть
неоднозначно определена.

{\bf Замечание.} В общем случае, последовательность наилучших
приближений может быть как конечной так и бесконечной. Но в данной
работе мы не будем сталкиваться с конечными последовательностями
наилучших приближений.

{\bf Замечание.} Функция $\psi_{\mathbf{v}}(t)$ есть кусочно
постоянная убывающая функция:
$$
\psi_{\mathbf{v}}(t)=\sqrt{\|q_n v_1\|^2+\|q_n v_2\|^2}=\sqrt{|q_n
v_1-p_{n,1}|^2+|q_n v_2-p_{n,2}|^2}=|q_n\mathbf{v}-\mathbf{p}_n|,
\mbox{ при } q_n\leqslant t<q_{n+1}
$$

В силу вышеприведенного замечания можно дать еще одно определение
двумерного спектра Дирихле, в евклидовой норме

{\bf Определение.}
$$
\mathbb{D}_2=\left\{\lambda\in\mathbb{R}\mid\exists\mathbf{v}\in\mathbb{R}^2:\;
\limsup\limits_{n\rightarrow\infty}
{q_{n+1}\cdot|q_n\mathbf{v}-\mathbf{p}_n|^2}=\lambda\right\}.
$$

{\bf Определение.} Пусть даны вектор ${\bf{v}}\in\mathbb{R}^2$,
$Q>0$ и $R>0$. Мы определим цилиндр $\Pi$:
$$
\Pi({\bf{v}}, Q, R)=
\left\{(q,{\bf{p}})\in\mathbb{R}\times\mathbb{R}^{2}\mid
q\in[0,Q],\; |q{\bf{v}} - {\bf{p}}|\leqslant R\right\}.
$$
Обозначим\\
через $\mbox{vol}(\Pi)=\pi Q R^2$ --- объем цилиндра $\Pi$,\\
через $\mbox{int}\Pi$ --- множество внутренних точек цилиндра
$\Pi$,\\
через $\partial\Pi$
--- границу цилиндра $\Pi$.\\
Положим
$$
\overline{\Pi}=
\left\{(q,{\bf{p}})\in\mathbb{R}\times\mathbb{R}^{2}\mid
q\in\mathbb{R},\; |q{\bf{v}} - {\bf{p}}|\leqslant R\right\}.
$$
Передней гранью цилиндра $\Pi$ будем называть множество
$$
\left\{(q,{\bf{p}})\in\mathbb{R}\times\mathbb{R}^{2}\mid q=Q,\;
|q{\bf{v}} - {\bf{p}}|< R\right\}.
$$
Направляющей цилиндра $\Pi$
будем называть вектор
$(1,{\bf{v}})\in\mathbb{R}\times\mathbb{R}^{2}$, а также любой
ненулевой вектор ему сонаправленный.\\
Боковой поверхностью цилиндра
$\Pi$ будем называть множество
$$\left\{(q,{\bf{p}})\in\mathbb{R}\times\mathbb{R}^{2}\mid
q\in(0,Q),\; |q{\bf{v}} - {\bf{p}}|= R\right\},
$$
а продолжением
боковой поверхности цилиндра $\Pi$ будем называть множество
$$
\left\{(q,{\bf{p}})\in\mathbb{R}\times\mathbb{R}^{2}\mid
q\in\mathbb{R},\; |q{\bf{v}} - {\bf{p}}|= R\right\}.
$$
Высотой или длиной цилиндра будем называть число $q$.\\
Радиусом цилиндра $\Pi$ будем называть число $R$.

{\bf Определение.} Для вектора $\mathbf{v}\in\mathbb{R}^2$ и
последовательности
$$
\mathcal{Z}:\mathbf{w}_n=(q_n,\mathbf{p}_n)\in\mathbb{N}\times\mathbb{Z}^{2},\qquad
n\geqslant0,
$$
мы положим $R_0=1$, а при $n\geqslant1$ определим
$R_n=|q_{n-1}\mathbf{v}-\mathbf{p}_{n-1}|$, $\Pi_n=\Pi(\mathbf{v},
q_n, R_n)$, $V_n=\mbox{vol}(\Pi_n)=\pi
q_{n}\cdot|q_{n-1}\mathbf{v}-\mathbf{p}_{n-1}|^2$.

{\bf Замечание.} Здесь и далее мы будем пользоваться введенными в
вышеприведенных определениях обозначениями.

В силу вышеприведенных обозначений можно дать еще одно определение
двумерного спектра Дирихле, в евклидовой норме

{\bf Определение.}
$$
\mathbb{D}_2=\left\{\lambda\in\mathbb{R}\mid\exists\mathbf{v}\in\mathbb{R}^2:\;
\limsup\limits_{n\rightarrow\infty}
\frac{1}{\pi}V_{n+1}=\lambda\right\}.
$$

Мы докажем теорему 2 из которой очевидно, в силу вышеприведенного
определения, следует теорема 1.

{\bf Теорема 2.}\\
Пусть $\lambda\in\left[0,\frac{2}{\sqrt{3}}\right]$. Тогда
существует континуальное множество векторов
$\mathbf{v}\in\mathbb{R}^2$ таких что
$$
\lim\limits_{n\rightarrow\infty}\frac{1}{\pi}V_{n+1}=\lambda
$$

Фактически в работе доказан более общий результат.

{\bf Теорема 3.}\\
Пусть $\{\Delta_n\}_{n=1}^{\infty}$ --- произвольная
последовательность отрезков из отрезка
$\left[0,\frac{2}{\sqrt{3}}\right]$.\\
Тогда существует континуальное
множество векторов $\mathbf{v}\in\mathbb{R}^2$ таких что
$$
\frac{1}{\pi}V_{n}\in\Delta_n\qquad \forall n\in\mathbb{N}
$$

{\bf Замечание.} Из теоремы 3 очевидно следует теорема 2.
Действительно, пусть $\lambda\in\left[0,\frac{2}{\sqrt{3}}\right]$.
Для доказательства теоремы 2 достаточно в теореме 3 определить
последовательность открытых интервалов следующим образом
$\Delta_n=\left[\lambda-\frac{1}{n},\lambda+\frac{1}{n}\right]\cap\left[0,\frac{2}{\sqrt{3}}\right]\qquad
\forall n\in\mathbb{N}$.

Теперь дадим еще одно определение последовательности наилучших
приближений.

{\bf Определение.} Последовательность векторов
$$
\mathcal{Z}:\mathbf{w}_n=(q_n,\mathbf{p}_n)\in\mathbb{Z}^{3},\qquad
n\geqslant0
$$
называется последовательностью наилучших приближений для вектора
$\mathbf{v}\in\mathbb{R}^2$, в евклидовой норме, если\\
1) $q_0=1$\\
2) $\left(\mbox{int}\Pi_n\right)\cap\mathbb{Z}^3=\varnothing$ $(n\geqslant 0)$\\
3) $q_{n+1}>q_n$ $(n\geqslant 0)$\\
4) $R_{n+1}<R_n$ $(n\geqslant 0)$

{\bf Замечание.} $\mathbf{w}_{n-1},\mathbf{w}_n\in\partial\Pi_n$, и
более того, точка $\mathbf{w}_{n-1}$ лежит на боковой поверхности
цилиндра $\Pi_n$, а точка $\mathbf{w}_n$ лежит в передней грани
цилиндра $\Pi_n$.

{\bf Доказательство теоремы 3.}

Мы построим последовательность векторов
$$
\mathcal{Z}:\mathbf{w}_n=(q_n,\mathbf{p}_n)\in\mathbb{Z}^{3},\qquad
n\geqslant0
$$
обладающих нижеследующими свойствами.

Для простоты изложения введем вспомогательные обозначения:\\
$\mathbf{v}_n=\frac{\mathbf{p}_n}{q_n}$, $q_{-1}=0$,
$\mathbf{p}_{-1}=(1,0)$,\\
при $0 \leqslant \nu \leqslant n$ положим
$R_n^\nu=\left|\frac{q_{\nu-1}}{q_n}\mathbf{p}_n-\mathbf{p}_{\nu-1}\right|$,
$\Pi_n^\nu=\Pi\left(\mathbf{v}_n, q_\nu,
R_n^\nu\right)$, $V_n^\nu=\mbox{vol}(\Pi_n^\nu)$. \\

 1) $\Pi_n^\nu\cap\mathbb{Z}^3=\varnothing$ $(0 \leqslant \nu \leqslant n)$

 2) $q_{n}>q_{n-1}$ $(n\geqslant 1)$

 3) $R_{n}^{\nu}<\frac{1}{2}R_n^{\nu-1}$ $(1 \leqslant \nu \leqslant n)$

 4) $\frac{1}{\pi}V_{n}^\nu\in\mbox{int}\Delta_{\nu}$ $(1 \leqslant \nu \leqslant n)$

 5) $\left|\mathbf{v}_n-\mathbf{v}_{n-1}\right|<\frac{1}{2^n}$ $(n\geqslant 1)$

 6) $|R_{n}^{\nu}-R_{n-1}^{\nu}|<\frac{1}{2^n}$ $(1 \leqslant \nu \leqslant n-1)$

Мы построим последовательность векторов $\mathcal{Z}$ индуктивным
образом (индукция по $n$).

1) Базис индукции: $n=0$

Определим вектор ${\bf{w}}_0$ так: $q_0=1$, ${\bf{p}}_0=(0,0)$.
Непосредственной проверкой легко убедиться, что при $n=0$ все
свойства 1) - 6) выполняются.

2) Предположение индукции: пусть уже имеется последовательность
векторов ${\bf{w}}_0,\dots,{\bf{w}}_{n-1}$, такая что свойства 1) -
6) выполнены.

3) Индуктивный шаг: построим вектор ${\bf{w}}_{n}$ и докажем, что
свойства 1) - 6) для него выполнены.

Пусть $\pi_{n-1}$ --- произвольная вполне рациональная
гиперплоскость в пространстве $\mathbb{R}^3$ проходящая через начало
координат и точку ${\bf{w}}_{n-1}$.

Пусть $\pi'_{n-1}$ --- соседняя к $\pi_{n-1}$ и параллельная ей
вполне рациональная гиперплоскость.

{\bf Замечание.} Соседних к $\pi_{n-1}$ гиперплоскостей две, в
качестве $\pi'_{n-1}$ мы можем взять любую из них. Это замечание
пригодится в дальнейшем для доказательства континуальности множества
векторов $\mathbf{v}\in\mathbb{R}^2$.

Пусть $G$ --- линейное отображение в $\mathbb{R}^3$, которое
однозначно задается следующими свойствами:

 1) $|G|=1$

 2) $G({\bf{w}}_{n-1})=(q_{n-1}, 0, 0)$

 3) $G(\pi_{n-1})=\{ z = 0 \}$

 4) $G(\sigma)=\sigma$, где $\sigma=\{x=1\}$

 5) $G(\pi'_{n-1})=\{ z = h \}$ и $h > 0$

Для простоты изложения данной части доказательства введем
сокращенные обозначения:
 $q = q_{n-1}$,
 $\widetilde{{\bf{w}}}_{n-1} = G({\bf{w}}_{n-1})$,
 $\widetilde{\pi}=\widetilde{\pi}_{n-1}= G(\pi_{n-1})$,
 $\widetilde{\pi}'=\widetilde{\pi}'_{n-1}= G(\pi'_{n-1})$,
 $\pi=\pi_{n-1}$,
 $\pi'=\pi'_{n-1}$,
 $\Lambda=\mathbb{Z}^3$,
 $\widetilde{\Lambda}=G(\Lambda)$,
 $\Gamma=\pi\cap\Lambda$,
 $\widetilde{\Gamma}=G(\Gamma)=\widetilde{\pi}\cap\widetilde{\Lambda}$

{\bf Свойства.}

 1) $\widetilde{\pi}$ и $\widetilde{\pi}'$ -- соседние вполне
 рациональные (относительно решетки $\Lambda$ ) гиперплоскости,
 а $h$~---~евклидово расстояние между ними.

 2) Решетка $\widetilde{\Gamma}$ --- двумерная и $\det\widetilde{\Gamma}=\frac{1}{h}$.
 Все точки решетки $\widetilde{\Gamma}$ располагаются на параллельных прямых
 $l_k=\{(x,y,z)\mid y=kd,z=0\}$ ($k\in\mathbb{Z}$) с шагом $q$.
 Расстояние между соседними прямыми равно $d$.

 3) Решетка $\widetilde{\Gamma}'$ --- двумерная и $\det\widetilde{\Gamma}'=\frac{1}{h}$.
 Все точки решетки $\widetilde{\Gamma}$ располагаются на параллельных прямых
 $l'_k=\{(x,y,z)\mid y=kd + b,z=h\}$ ($k\in\mathbb{Z}$, а $b$ --- некоторое действительное число) с шагом $q$.
 Расстояние между соседними прямыми равно $d$.

 4) $\det\widetilde{\Lambda}=qhd=1$.

 5) Преобразование $G$ не изменяет координату по оси $x$.

 6) Преобразование $G$ переводит любой цилиндр $\Pi$ в цилиндр $\widetilde{\Pi}=G(\Pi)$.
 При этом сохраняются радиус, длина и объем цилиндра $\Pi$.
 Будем пользоваться этим обозначением в этой части доказательства.

 7) Преобразование $G$ не только переводит плоскость $\sigma$ в себя $G(\sigma)=\sigma$,
 но и сохраняет евклидовы расстояния между точками в плоскости $\sigma$.

Переформулируем задачу построение вектора
${\bf{w}}_{n}=(q_n,\mathbf{p}_n)\in\mathbb{N}$:

Требуется построить последовательность векторов
$$
\widetilde{{\bf{w}}}_0=(q_0,\widetilde{\mathbf{p}}_0), \dots,
\widetilde{{\bf{w}}}_{n}=(q_{n},\widetilde{\mathbf{p}}_{n})\in\widetilde{\Lambda},
$$
таких что нижеприведенные свойства выполнены.

Для простоты изложения введем вспомогательные обозначения:\\
$\widetilde{\mathbf{v}}_n=\frac{\widetilde{\mathbf{p}}_n}{q_n}$,
$q_{-1}=0$,
$\mathbf{p}_{-1}=(1,0)$,\\
при $0 \leqslant \nu \leqslant n$ положим
$R_n^\nu=\left|\frac{q_{\nu-1}}{q_n}\mathbf{p}_n-\mathbf{p}_{\nu-1}\right|
=\left|\frac{q_{\nu-1}}{q_n}\widetilde{\mathbf{p}}_n-\widetilde{\mathbf{p}}_{\nu-1}\right|$,
$\widetilde{\Pi}_n^\nu=\Pi\left(\widetilde{\mathbf{v}}_n, q_\nu,
R_n^\nu\right)$

 1) $\widetilde{\Pi}_n^\nu\cap\widetilde{\Lambda}=\varnothing$ $(0 \leqslant \nu \leqslant n)$

 2) $q_{n}>q_{n-1}$ $(n\geqslant 1)$

 3) $R_{n}^{\nu}<\frac{1}{2}R_n^{\nu-1}$ $(1 \leqslant \nu \leqslant n)$

 4) $\frac{1}{\pi}V_{n}^\nu\in\mbox{int}\Delta_{\nu}$ $(1 \leqslant \nu \leqslant n)$

 5) $\left|\mathbf{v}_n-\mathbf{v}_{n-1}\right|<\frac{1}{2^n}$ $(n\geqslant 1)$

 6) $|R_{n}^{\nu}-R_{n+1}^{\nu}|<\frac{1}{2^n}$ $(1 \leqslant \nu\leqslant n)$

Причем вектора
$$
\widetilde{{\bf{w}}}_0=(q_0,\widetilde{\mathbf{p}}_0)=G({{\bf{w}}}_0),
\dots,
\widetilde{{\bf{w}}}_{n-1}=(q_{n-1},\widetilde{\mathbf{p}}_{n-1})=G({{\bf{w}}}_{n-1}),
$$
уже построены и для них, в силу свойств отображения $G$, уже
выполнены свойства 1) - 6).

Необходимо построить вектор
$\widetilde{{\bf{w}}}_{n}=(q_n,\widetilde{\mathbf{p}}_n)\in\widetilde{\Lambda}$,
такой что выполнены свойства 1) - 6).

Как только мы найдем такой вектор, то положим
${\bf{w}}_{n}=G^{-1}(\widetilde{{\bf{w}}}_{n})$ и тем самым
построение последовательности $\mathcal{Z}$ будет завершено.

Для дальнейшего доказательства нам понадобится определить три
вспомогательных множества $A_0, A_1, A_2\subset\mathbb{R}^3$.

{\bf Замечание.} Элементы множеств $A_0, A_1$ и $A_2$ мы будем, для
простоты изложения, обозначать через $(x_0, y_0, z_0)$, $(x_1, y_1,
z_1)$ и $(x_2, y_2, z_2)$ соответственно.

Определим множество
$$
A_0 = \{ (x_0, y_0, z_0) \in \mathbb{R}^3 \mid x_0=q, y_0>0, z_0>0
\}.
$$

Каждой точке множества $A_0$ поставим в соответствие цилиндр
$$
\Pi(0)(x_0, y_0,
z_0)=\Pi\left(\left(\frac{y_0}{q},\frac{z_0}{q}\right),
\frac{qh}{z_0}, \sqrt{y_0^2+z_0^2}\right).
$$

{\bf Замечание.} Нижеприведенные свойства цилиндра $\Pi(0)(x_0, y_0,
z_0)$ однозначно его задают.

1) Вектор $(x_0, y_0, z_0)$ является направляющим для цилиндра
$\Pi$.

2) Центр передней грани цилиндра, лежит в плоскости
$\widetilde{\pi}'$.

3) Точка $\widetilde{{\bf{w}}}_{n-1} = (q, 0, 0)$ лежит на боковой
поверхности цилиндра $\Pi$ (или на ее продолжении).

Определим множество
$$
A_1 = \{ (x_1, y_1, z_1) \in \mathbb{R}^3 \mid x_1>0, y_1>0, z_1=h
\}.
$$

Каждой точке множества $A_1$ поставим в соответствие цилиндр
$$
\Pi(1)(x_1, y_1,
z_1)=\Pi\left(\left(\frac{y_1}{x_1},\frac{h}{x_1}\right), x_1,
\frac{q}{x_1}\sqrt{y_1^2+h^2}\right).
$$

{\bf Замечание.} Нижеприведенные свойства цилиндра $\Pi(1)(x_1, y_1,
z_1)$ однозначно его задают.

1) Вектор $(x_1, y_1, z_1)$ является направляющим для цилиндра
$\Pi$.

2) Центр передней грани цилиндра, лежит в плоскости
$\widetilde{\pi}'$.

3) Точка $\widetilde{{\bf{w}}}_{n-1} = (q, 0, 0)$ лежит на боковой
поверхности цилиндра $\Pi$ (или на ее продолжении).

Определим множество
$$
A_2 = \{ (x_2, y_2, z_2) \in \mathbb{R}^3 \mid x_2>0, y_2>0, z_2=0
\}.
$$

Каждой точке множества $A_2$ поставим в соответствие цилиндр
$$
\Pi(2)(x_2, y_2, z_2)=\Pi\left(\left(\frac{x_2 y_2}{x_2^2 +
q^2},\frac{q y_2}{x_2^2 + q^2}\right), \frac{h(x_2^2 + q^2)}{q y_2},
\frac{q y_2}{\sqrt{x_2^2 + q^2}} \right).
$$

{\bf Замечание.} Нижеприведенные свойства цилиндра $\Pi(2)(x_2, y_2,
z_2)$ однозначно его задают.

1) Прямая $\{x=0, y=y_2, z=0\}$ касается боковой поверхности
цилиндра $\Pi$ (или ее продолжения) в точке $(x_2, y_2, z_2)$.

2) Центр передней грани цилиндра, лежит в плоскости
$\widetilde{\pi}'$.

3) Точка $\widetilde{{\bf{w}}}_{n-1} = (q, 0, 0)$ лежит на боковой
поверхности цилиндра $\Pi$ (или на ее продолжении).

Зададим тройку согласованных биективных отображений между
множествами $A_0, A_1$ и $A_2$:
$$
\left\{
\begin{array}{l}
x_0 = q \\
y_0 = \frac{q y_1}{x_1} \\
z_0 = \frac{q h}{x_1}
\end{array}
\right.
\left\{
\begin{array}{l}
x_1 = \frac{q h}{z_0} \\
y_1 = \frac{q y_0}{z_0} \\
z_1 = h
\end{array}
\right.
\left\{
\begin{array}{l}
x_1 = \frac{h (x_2^2 + q^2)}{y_2 q} \\
y_1 = \frac{h x_2}{q} \\
z_1 = h
\end{array}
\right.
\left\{
\begin{array}{l}
x_2 = \frac{q y_1}{h} \\
y_2 = \frac{q (y_1^2 + h^2)}{x_1 h} \\
z_2 = 0
\end{array}
\right.
\left\{
\begin{array}{l}
x_0 = q \\
y_0 = \frac{q x_2 y_2}{x_2^2 + q^2} \\
z_0 = \frac{q^2 y_2}{x_2^2 + q^2}
\end{array}
\right.
\left\{
\begin{array}{l}
x_2 = \frac{q y_0}{z_0} \\
y_2 = \frac{y_0^2 + z_0^2}{z_0} \\
z_2 = 0
\end{array}
\right.
$$

{\bf Замечание.} Легко проверить, что заданные преобразования
действительно определяют тройку согласованных биективных отображений
между множествами $A_0, A_1$ и $A_2$.

{\bf Замечание.} Так же легко проверить и то, что при заданных
преобразованиях инвариантом является цилиндр соответствующий
элементам множеств $A_0, A_1$ и $A_2$:
$$
\Pi(0)(x_0, y_0, z_0)=\Pi(1)(x_1, y_1, z_1)=\Pi(2)(x_2, y_2, z_2).
$$

{\bf Замечание.} В дальнейшем для простоты изложения мы будем
отождествлять соответствующие элементы множеств $A_0, A_1$ и $A_2$.
Будем считать, что множества $A_0, A_1$ и $A_2$ задают различные
параметризации одного и того же семейства цилиндров.

Легко проверить, что верна следующая лемма.

{\bf Лемма.} Пусть $r>0$. Тогда следующие уравнения равносильны

1) $x_1 = \frac{q (y_1^2 + h^2)}{2 r h}$

2) $y_2 = 2r$

3) $y_0^2 + (z_0 - r)^2 = r^2$

и задают множество цилиндров с объемом $V=V(r)=2\pi rqh$

{\bf Замечание.} Уравнение $V=V(r)=2\pi rqh$ можно переписать в виде
$V(r)=\frac{2\pi r}{d}$

{\bf Замечание.} Уравнение $y_0^2 + (z_0 - r)^2 = r^2$ можно
переписать в виде $\frac{y_0^2+z_0^2}{z_0} = 2r$

Вы берем параметр $r$, таким образом, чтобы
$\lambda^*=\frac{2r}{d}\in\Delta_n$ и чтобы число $\frac{2rh}{d^2}$
было иррациональным.

Определим множество
$$
B_2 = \{ (x_2, y_2, z_2) \in \mathbb{R}^3 \mid x_2>0, y_2>0, z_2=0,
\overline{\Pi}(2)(x_2, y_2, z_2)\cap\widetilde{\Gamma}=\varnothing
\}.
$$

{\bf Замечание.} Если $(x_2, y_2, z_2)\in B_2$, то $\Pi(2)(x_2, y_2,
z_2)\cap\Gamma=\varnothing$

{\bf Замечание.} $B_2$ --- открытое множество.

{\bf Замечание.} Если $(x_2, y_2, z_2)\in B_2\cup\partial B_2$, то
$\theta(x_2, y_2, z_2)\in B_2\;\;\;\forall\theta\in(0,1)$

{\bf Замечание.} Если $(x_2, y_2, z_2)\in B_2$, то
$(x_2+k\frac{y_2q}{d}, y_2, z_2)\in B_2\;\;\;\forall k\geqslant0$

{\bf Замечание.} Существует такое $0<a\leqslant q$, что $(a, d,
0)\in\Gamma$.

{\bf Лемма.} Существует $\varepsilon>0$, такое что для любого
$k\geqslant1$ множество $B_2$ содержит в себе прямоугольник
ограниченный заданными прямыми:
$$
y_2=\lambda^*d
$$
$$
y_2=\lambda^*d-\frac{\lambda^*d\varepsilon}{\lambda^*\left(a +
\left(k+\frac{1}{2}\right)q\right)+\varepsilon}
$$
$$
x_2=\lambda^*\left(a + \left(k+\frac{1}{2}\right)q\right)
$$
$$
x_2=\lambda^*\left(a +
\left(k+\frac{1}{2}\right)q\right)-\varepsilon
$$

{\bf Доказательство леммы.} Легко проверить, что
$$
\frac{2}{\sqrt{3}}\frac{(a, d, 0)+(a + q, d,
0)}{2}=\frac{2}{\sqrt{3}}\left(a + \frac{1}{2}q, d, 0\right)\in
\partial B_2,
$$
откуда следует, что
$$
\lambda^*\left(a + \frac{1}{2}q, d, 0\right)\in B_2
$$
Так как множество $B_2$ открыто, то
$$
\exists
\varepsilon>0\;\;\;\forall\theta\in[-\varepsilon,\varepsilon]\;\;\;
\lambda^*\left(a + \frac{1}{2}q+\theta, d, 0\right)\in B_2.
$$
В силу свойств множества $B_2$ получаем
$$
\exists \varepsilon>0\;\;\;\forall
k\geqslant0\;\;\;\forall\theta\in[-\varepsilon,\varepsilon]\;\;\;
\lambda^*\left(a + \left(k+\frac{1}{2}\right)q+\theta, d,
0\right)\in B_2.
$$
Отсюда получаем, что треугольник с вершинами в
точках\\
$(0,0,0)$, $\lambda^*\left(a +
\left(k+\frac{1}{2}\right)q-\varepsilon, d, 0\right)$ и
$\lambda^*\left(a + \left(k+\frac{1}{2}\right)q+\varepsilon, d,
0\right)$ содержится (за исключением вершины $(0,0,0)$) в множестве
$B_2$.

Легко проверить, что вышеупомянутый треугольник содержит в себе
прямоугольник ограниченный заданными прямыми:
$$
y_2=\lambda^*d
$$
$$
y_2=\lambda^*d-\frac{\lambda^*d\varepsilon}{\lambda^*\left(a +
\left(k+\frac{1}{2}\right)q\right)+\varepsilon}
$$
$$
x_2=\lambda^*\left(a + \left(k+\frac{1}{2}\right)q\right)
$$
$$
x_2=\lambda^*\left(a +
\left(k+\frac{1}{2}\right)q\right)-\varepsilon
$$
Что и требовалось доказать. Лемма доказана.

Рассмотрим семейство полученных в лемме прямоугольников в терминах
множества $A_1$

$$
x_1=\frac{q (y_1^2 + h^2)}{\lambda^* dh}
$$
$$
x_1=\frac{q (y_1^2 +
h^2)}{\left(\lambda^*d-\frac{\lambda^*d\varepsilon}{\lambda^*\left(a
+ \left(k+\frac{1}{2}\right)q\right)+\varepsilon}\right) h}
$$
$$
y_1=\lambda^*\left(a + \left(k+\frac{1}{2}\right)q\right)\frac{h}{q}
$$
$$
y_1=\lambda^*\left(a +
\left(k+\frac{1}{2}\right)q\right)\frac{h}{q}-\varepsilon\frac{h}{q}
$$

Первые два уравнения задают непересекающиеся параболы с общей осью,
и при выборе достаточно большого параметра $k$ расстояния между
ветвями парабол вдоль оси $x$ становятся сколь угодно большими (и в
частности, больше числа $q$) в заданном диапазоне по оси $y$:
$$
\frac{q (y_1^2 +
h^2)}{\left(\lambda^*d-\frac{\lambda^*d\varepsilon}{\lambda^*\left(a
+ \left(k+\frac{1}{2}\right)q\right)+\varepsilon}\right) h}-\frac{q
(y_1^2 + h^2)}{\lambda^*d h}=
$$
$$
=\frac{q (y_1^2 + h^2)}{\lambda^*d h}
\left(\frac{1}{\left(1-\frac{\varepsilon}{\lambda^*\left(a +
\left(k+\frac{1}{2}\right)q\right)+\varepsilon}\right)}-1\right)\rightarrow+\infty\;\;\;(k\rightarrow+\infty)
$$

Оставшиеся два уравнения задают пару параллельных прямых идущих
вдоль оси $x$, расстояние между которыми равно константе
$\varepsilon\frac{h}{q}$ не зависящей от выбора параметра $k$.

Итак, полученные в лемме прямоугольники, в терминах множества $A_1$
представляют собой периодически идущие куски полосок постоянной
ширины (с периодом $\lambda^* h=\frac{2rh}{d}$ по оси $y$)
расположенные в плоскости $\widetilde{\pi}'$ вдоль оси $x$,
ограниченных парой парабол.

В силу свойств отображения $G$, в плоскости $\widetilde{\pi}'$ все
точки решетки $\widetilde{\Gamma}$ располагаются на прямых
(периодично с шагом $q$), параллельных оси $x$, и расстояние между
соседними прямыми равно~$d$.

В силу выбора параметра $r$ отношение периодов $\frac{2rh}{d}$ и $d$
иррационально, и значит существует сколь угодно много возможностей
выбрать параметр $k$ так, чтобы соответствующий кусок полоски в
множестве $A_1$ содержал точку (или несколько точек) решетки
$\widetilde{\Gamma}$. Параметр $k$ в этом случае будем называть
допустимым.

Выберем одну из соответствующих допустимому параметру $k$ точек
решетки $\widetilde{\Gamma}$ и сопоставим ее в соответствие этому
параметру $k$.

В дальнейшем будем говорить, что достаточно большому подходящему
параметру $k$ соответствует точка решетки $\widetilde{\Gamma}$, или
что достаточно большому подходящему параметру $k$ соответствует
параметр из множества $A_0$, $A_1$ или $A_2$.

Заметим, что при выборе достаточно большого допустимого параметра
$k$ аргумент $x_1$ соответствующего параметра $(x_1, y_1, z_1)\in
A_1$ неограниченно возрастает.

Рассмотрим теперь семейство полученных в лемме прямоугольников в
терминах множества~$A_0$
$$
y_0^2 + \left(z_0 - \frac{1}{2}\lambda^*d\right)^2 =
\left(\frac{1}{2}\lambda^*d\right)^2
$$
$$
y_0^2 + \left(z_0 -
\frac{1}{2}\left(\lambda^*d-\frac{\lambda^*d\varepsilon}{\lambda^*\left(a
+ \left(k+\frac{1}{2}\right)q\right)+\varepsilon}\right)\right)^2 =
\left(\frac{1}{2}\left(\lambda^*d-\frac{\lambda^*d\varepsilon}{\lambda^*\left(a
+ \left(k+\frac{1}{2}\right)q\right)+\varepsilon}\right)\right)^2
$$
$$
y_0=\frac{z_0}{q}\left(\lambda^*\left(a +
\left(k+\frac{1}{2}\right)q\right)\right)
$$
$$
y_0=\frac{z_0}{q}\left(\lambda^*\left(a +
\left(k+\frac{1}{2}\right)q\right)-\varepsilon\right)
$$
Полученные в лемме прямоугольники, в терминах множества $A_0$
представляют собой области в плоскости $\{x=q\}$, ограниченные двумя
окружностями, которые касаются друг друга внутренним образом в точке
${\bf{w}}_{n-1}$ и двумя лучами исходящими из начала координат.

Легко проверить, что при выборе достаточно большого допустимого
параметра $k$, соответствующий параметр $(x_0, y_0, z_0)\in A_0$,
находится в сколь угодно малой окрестности точки ${\bf{w}}_{n-1}$.

Докажем, что можно выбрать такой достаточно большой допустимый
параметр $k$ с соответствующим вектором (точкой решетки)
$\widetilde{{\bf{w}}}_{n}$ так чтобы свойства 1) - 6) были
выполнены:

По построению свойство 4) выполнено.

Так как по построению
$\widetilde{\Pi}_n^n\cap\widetilde{\Gamma}=\varnothing$, то в силу
симметричного расположение цилиндра $\widetilde{\Pi}_n^n$
относительно решетки $\widetilde{\Lambda}$, следует что
$\widetilde{\Pi}_n^n\cap\widetilde{\Gamma}'=\varnothing$. Так как
при выборе достаточно большого допустимого параметра $k$
соответствующий параметр $(x_0, y_0, z_0)\in A_0$ находится в сколь
угодно малой окрестности точки ${\bf{w}}_{n-1}$, то можно добиться
выполнения свойств 1) (при $0 \leqslant \nu < n$), 3), 5), 6) и того
чтобы радиус $R_n^n$ был сколь угодно малым, и как следствие
$\widetilde{\Pi}_n^n\cap\widetilde{\Lambda}=\varnothing$, то есть
выполнения свойства 1).

Так как при выборе достаточно большого допустимого параметра $k$
аргумент $x_1$ соответствующего параметра $(x_1, y_1, z_1)\in A_1$
неограниченно возрастает, то можно добиться выполнения свойства 2).

Индуктивное построение завершено.

Доказательство леммы завершено.

В силу свойств построенной последовательности $\mathcal{Z}$ и
критерия Коши сходимости последовательности получаем, что
$\exists\lim\limits_{i\rightarrow\infty}{\bf{v}}_{i}={\bf{v}}$.

Докажем, что построенная последовательность векторов
$$
\mathcal{Z}:\mathbf{w}_n=(q_n,\mathbf{p}_n)\in\mathbb{Z}^{3},\qquad
n\geqslant0
$$
является последовательностью наилучших приближения для вектора
$\mathbf{v}\in\mathbb{R}^2$.

Действительно,

 1) $q_0=1$ -- это выполнено по построению

 2) $\left(\mbox{int}\Pi_n\right)\cap\mathbb{Z}^3=\varnothing$
 $(n\geqslant 0)$ -- это выполнено по построению, так как
 $\Pi_n=\lim\limits_{\nu\rightarrow\infty}\Pi_n^\nu$

 3) $q_{n+1}>q_n$ $(n\geqslant 0)$ -- это выполнено по построению

 4) $R_{n+1}<R_n$ $(n\geqslant 0)$ -- это выполнено по построению

Континуальность множества векторов $\mathbf{v}\in\mathbb{R}^2$
обеспечивается следующим образом. На каждом шаге индуктивного
процесса основной леммы, при применении леммы 2 проводим построение
не только для одной соседней гиперплоскости $\pi'_{n-1}$, но и для
второй. Нетрудно показать, что получающиеся при этом векторы
$\mathbf{v}\in\mathbb{R}^2$ будут различны.

Доказательство теоремы 3 завершено.

\end{document}